\documentclass[11pt]{article}   
\usepackage{amsfonts, amsthm, amsmath}   
\theoremstyle{plain}
\newtheorem{theorem}{Theorem}[section]

\theoremstyle{definition}
\newtheorem{definition}[theorem]{Definition}

\newtheorem{problem}[theorem]{Problem}
\newtheorem{fromnowon}[theorem]{From now on}

\theoremstyle{remark}

\newenvironment{lisp}{\begingroup\footnotesize}{\endgroup}

\newcommand{\boxtt}[1]{\mbox{\texttt{#1}}}
\newcommand{\iten}{\item\vspace{-6pt}}
\newcommand{\N}{\mbox{$\mathbb{N}$}}
\newcommand{\R}{\mbox{$\mathbb{R}$}}
\newcommand{\Star}{{\textstyle \ast}}

\newcommand{\Z}{\mbox{$\mathbb{Z}$}}

\title{Constructive Algebraic Topology\footnote{This text was used as a
background paper for a plenary talk of the second author during the EACA
Congress of Tenerife, September 1999. A ``general public'' version  has
appeared in~\cite{SRGR5}; it is an excellent introduction for the present
text.}}
\author{\small \em Julio Rubio, Francis Sergeraert }
\date{}

\begin{document}
\maketitle

\begin{abstract} The classical ``computation'' methods in Algebraic
Topology most often work by means of highly infinite objects and in fact
\emph{are not} constructive. Typical examples are shown to describe the nature
of the problem. The Rubio-Sergeraert solution for Constructive Algebraic
Topology is recalled. This is not only a theoretical solution: the concrete
computer program \emph{Kenzo} has been written down which precisely follows
this method. This program has been used in various cases, opening new research
subjects and producing in several cases significant results unreachable by
hand. In particular the Kenzo program can compute the first homotopy groups of
a simply connected \emph{arbitrary} simplicial set. \end{abstract}

\section*{Introduction}

The computation of \emph{homotopy groups} in Algebraic Topology is known as a
difficult problem. Every pointed topological space $(X,x_0)$ has a family of
homotopy groups $\{\pi_n(X,x_0)\}_{n \geq 1}$, $\{\pi_n X\}$ in short, and
these groups are abelian for $n \geq 2$. The definition was given by Hurewicz
in 1935, and for the first non trivial space from this point of view, namely
the 2-sphere $S^2$, only the groups $\pi_2$ and $\pi_3$ were known at this
time, thanks to Hopf. The group $\pi_4 S^2 = \Z_2$ was determined by
Freudenthal in 1937. Thirteen years then passed without any new homotopy group
of sphere. The following groups $\pi_n S^2$ were obtained by Serre for $5 \leq
n \leq 9$, in 1950. In fact, for $n = 6$, Serre proved the group $\pi_6 S^2$
has twelve elements but did not succeed in choosing between both possible
solutions $\Z_{12}$ and $\Z_2 \oplus \Z_6$. Two years later, Barratt and
Paechter proved there exists an element of order 4 in $\pi_6 S^2$, so that
finally $\pi_6 S^2 = \Z_{12}$. See \cite[vol. I, pp 110 and 113]{SERR3} for
details and references.

More generally, Serre obtained a general \emph{finiteness result}.

\begin{theorem}
\emph{(Serre, \cite[p. 14 and pp. 171-207]{SERR3})} --- If $X$ is a simply
connected space such that the homology groups $H_n(X;\Z)$ are of finite type,
then the homotopy groups $\pi_n X$ are also abelian groups of finite type.
\end{theorem}

In particular the homology groups of simply connected finite polyhedra, for
example the simply connected compact manifolds, are of finite type, so that
their homotopy groups are also of finite type. Various methods allow to
combinatorially describe the finite polyhedra: these objects may be the
\emph{input} of an algorithm. An abelian group of finite type can also be
described by some character string: such a group could be the \emph{output} of
an algorithm. The following problem therefore makes sense.

\begin{problem}
Does there exist a general algorithm:
\begin{itemize}
\item
\textbf{Input:} A simply connected polyhedron $X$ and an integer $n \geq
2$;
\iten
\textbf{Output:} The homotopy group $\pi_n X$.
\end{itemize}
\end{problem}

A solution for this computability problem was given by Edgar Brown in 1956
\cite{BRWNE1}. He used the general organization just defined by Postnikov, now
known as the \emph{Postnikov tower}; then the result is not difficult when the
homology groups of the space $X$ are \emph{finite}; really finite, not only of
finite type: for example this simple method does not work for the 2-sphere
$S^2$ because the homology group $H_2 S^2 = \Z$ is of finite type (one
generator), but unfortunately is infinite. The difficult part of the work of
Edgar Brown then consisted in overcoming the birth of infinite objects in the
Postnikov tower. A complicated and tricky process was used to approximate these
infinite objects by finite ones and in this way Edgar Brown succeeded in
transforming the finiteness result of Serre into a computability result.

But let us quote Edgar Brown himself in the introduction of his article:

\begin{quotation}
\emph{It must be emphasized that although the procedures developed for solving
these problems are finite, they are much too complicated to be considered
practical.}
\end{quotation}

Forty years later this appreciation still holds, and will always hold, even
with the most powerful computer you can imagine: it is a consequence of the
hyper-exponential complexity of the algorithm designed by Edgar Brown.

The problem of finding new \emph{general} algorithms which on the contrary
could be \emph{concretely} used in significant cases was not seriously studied
up to 1985. This is so true that topologists from time to time meet some
difficulty in expressing precisely where the actual nature of a problem is,
when in fact it is a matter of computability. Section \ref{sec1} shows three
typical examples of this sort. Let us quote immediately another example found
in the introduction of \cite{GLMN}:
\begin{quotation}
\emph{The book by Cartan and Eilenberg contains essentially all the
constructions of homological algebra that constitute its computational tools,
namely standard resolutions and spectral sequences.}
\end{quotation}

Strictly speaking, this statement is correct, but it is also very misleading.
In the ``general'' domain of Homological Algebra, it is true, but if you intend
to apply these ``computational'' tools in Algebraic Topology, then you realize
an enormous \emph{gap} is in front of you, mainly when you have to determine
the higher differentials of the spectral sequences you are working with; and if
you succeed in finding them, a collection of hard extension problems can be
waiting at the abutment. The present paper is essentially devoted to these
questions.

One of the examples of Section \ref{sec1} asserts a computability problem in
homotopy theory is ``widely'' open. In fact three complete solutions are
available for several years. Section \ref{sec2} is devoted to a quick
description of these solutions, to their nature and what can be hoped about
their concrete use for computer calculations.

So far, only the Rubio-S. solution has led to a reasonably complete computer
program which has been used in significant cases. The main tool is standard
algebraic topology combined with \emph{functional programming} and Section
\ref{sec3} uses a didactic example to explain how a functional programming
method can be used to obtain efficient algorithms, even for solving problems
where there is a function neither in the input nor in the output.

The main ingredient in the Rubio-S. solution is the notion of \emph{object with
effective homology}. Such an object is a subtle combination via chain
equivalences of traditional \emph{effective} objects on one hand, and of other
\emph{locally effective} objects on the other hand. Section \ref{sec4}
describes the essential properties of these objects and what an object with
effective homology is.

The main tools of basic algebraic topology, mainly the Serre and
Eilenberg-Moore spectral sequences, may then be rewritten in such a way they
become \emph{algorithms} computing the desired homology groups when the
necessary data are given; such a property does not hold for the classical
spectral sequences. Section \ref{sec5} contains the main statements and
describes how they can be used for example to compute the homotopy groups of
simply connected simplicial sets with effective homology. A simple solution is
so obtained for the computability problem of homotopy groups; furthermore its
scope is much larger than Edgar Brown's one.

Section \ref{sec6} describes how these theoretical results led to a concrete
program named \emph{Kenzo}\footnote{Kenzo is the name of the author's
\emph{cat} and C.A.T. = Constructive Algebraic Topology; the next version of
our program will therefore be called \emph{Simba}, the daughter of Kenzo.}. It
is a Lisp program of 16000 lines (joint work with Xavier Dousson), now
www-available \cite{DSSS} with a rich documentation (340 pp) written by Yvon
Siret.

These results open new research fields; in Computer Science because of the
original type of functional programming which is required, but in theoretical
Algebraic Topology as well: the objects that are processed by the Kenzo program
are much too complicated to be studied by hand, specially around
\emph{algebraic fibrations}. These questions are considered in Section
\ref{sec7}.

Finally Section \ref{sec8} gives a few examples of calculations.

\section{Three examples.}\label{sec1}

A preprint by Karoubi \cite{KARB}, distributed in 1993, begins as follows:

\begin{quotation}\label{Karoubi}
\emph{The problem of finding a ``computable algebraic model'' for the homotopy
type of a CW-complex $X$ remains a widely open problem in topology.}
\end{quotation}

The notion of \emph{computable algebraic model} for a homotopy type is not
precisely defined in the text, but taking account of the rest of the paper, and
also of other related papers by the same author, it is clear the following
meaning is the right one:

\begin{definition}\label{computable-model}
A \emph{computable algebraic model} for the homotopy type of a space $X$ is an
additional structure $\mathcal{H}$ over the chain complex $C_\Star X$ such that
the pair $(C_\Star X, \mathcal{H})$ ``contains'' the homotopy type of $X$.
\end{definition}

Two spaces $X$ and $Y$ have the same homotopy type if there exist two
continuous maps $f:X \rightarrow Y$ and $g:Y \rightarrow X$ such that $g \circ
f$ and $f \circ g$ are homotopic to identity maps; from the point of view of
Algebraic Topology, both spaces are ``equal'', even if they are quite
different: for example a point and the infinite unit sphere $S^\infty \subset
\ell^2$ have the same homotopy type: this sphere is in fact contractible.

As usual, the additional data $\mathcal{H}$ must be \emph{natural} with respect
to $X$, that is, the mapping $X \mapsto (C_\Star X, \mathcal{H})$ should be a
functor. Several contexts are possible. If the chain complex $C_\Star X$ is the
\emph{singular} chain complex, then it is easy to give the required additional
structure (the canonical distinguished generators, namely the singular
simplices, and the simplicial operators), but the singular chain complex is a
functional space which is so enormous that no program can handle it: the object
so obtained is not \emph{computable}. The same in the simplicial context as
soon as the simplicial model is infinite, which is frequent. Karoubi wants a
chain complex of finite type in any dimension, for example the cellular chain
complex $C_\Star^{\mbox{\scriptsize cell}}(X)$ if~$X$ is presented as a
CW-complex of finite type in any dimension; much information about~$X$ is lost
in this chain complex and Karoubi searches an additional structure over this
chain complex which captures the homotopy type of $X$ at least. The structures
studied by Karoubi intensively use the notion of \emph{non-commutative
differential forms} and are interesting, but to our knowledge, the goal defined
by Karoubi is not yet reached by his method.

In fact three solutions now exist for Constructive Algebraic Topology, and two
of them exactly have the form that Karoubi looked for. In Justin Smith'
solution \cite{SMTH1,SMTH2}, the cellular chain complex
\(C_\Star^{\mbox{cell}}(X)\) is provided with a \(m\)-structure which, in
appropriate context, is a computable algebraic model for the homotopy type of
\(X\). In the Rubio-S. solution, the same chain complex is completed with two
other chain complexes and a few operators which give the same result. The
solution by Rolf Sch\"on \cite{SCHN} is not presented in this way but finally
is equivalent to both previous ones.

Let us quote now a paper by Carlsson and Milgram in James' Handbook of
Algebraic Topology \cite[p. 545]{CRML}:

\begin{quotation}
\emph{In Section 5 we showed that for a connected CW complex with no one cells
one may produce a CW complex, with cell complex given as the free monoid on
generating cells, each in one dimension less than the corresponding cell of
$X$, which is homotopy equivalent to [the loop space of $X$] $\Omega X$. To go
further one should study similar models for double loop spaces, and more
generally for iterated loop spaces.}

\emph{In principle this is direct. Assume $X$ has no $i$-cells for $1 \leq i
\leq n$ then we can iterate the Adams-Hilton construction of Section 5 and
obtain a cell complex which represents $\Omega^n X$. However the question of
determining the boundaries of the cells is very difficult as we already saw
with Adam's solution of the problem in the special case that $X$ is a
simplicial complex with $sk_1(X)$ collapsed to a point. It is possible to
extend Adams' analysis to $\Omega^2 X$, but as we will see there will be severe
difficulties with extending it to higher loop spaces except in the case where
$X = \Sigma^n Y$.}
\end{quotation}

The paper by Carlsson and Milgram is an excellent presentation of Adams' model
for a loop space of a simply connected CW complex and related questions. You
see the authors here consider a problem whose solution \emph{in principle is
direct}, but new \emph{severe} difficulties are soon announced which can in
fact be overcome only if the space $X$ in an iterate suspension $\Sigma^n X$.

In fact the actual problem is a \emph{computability} problem. The following
theorem can easily be deduced from Adams' construction. In the statement, the
operator $s^{-1}$ is the desuspension of the ``augmentation ideal'': the base
generator is removed and the degree $n$ of a generator becomes $n-1$; the
operator $T$ associates to a chain complex its tensor algebra, another chain
complex provided with a multiplicative structure.

\begin{theorem}
If $X$ is a CW complex with one 0-cell, without any $i$-cell ($1 \leq i \leq
n$), then \emph{there exists} for the chain complex: $$ G^n X = (Ts^{-1})^n
C^{\mbox{\scriptsize \emph{cell}}}_\Star(X) $$ a new differential $\delta$ such
that the chain complex $(G^n X, \delta)$ is the cellular chain complex of a CW
model of the iterate loop space $\Omega^n X$.
\end{theorem}

The \emph{existence} of the differential $\delta$ can be easily proved thanks
to Adams' work about the CW model of the first loop space (cf also
\cite{BAUS1}), but the existence proof is not constructive: it is made of a
mixture of combinatorial and topological arguments and certainly there are at
least ``severe difficulties'' to translate the topological constructions into
the combinatorial constructions that are necessary if you intend to obtain a
constructive existence proof for the differential $\delta$. The problem of
\emph{iterating the cobar construction} is the \emph{heart} of Algebraic
Topology: the main computability problems can be reduced to this one, and it is
not amazing this problem is a little severe. The three current solutions
\cite{SCHN, SRGR3,DSSS,SMTH1,SMTH2} for Constructive Algebraic Topology are
firstly solutions for the problem of iterating the cobar construction.

John McCleary tries in his book \cite{MCCL} to express the same idea in the
context of spectral sequences:

\begin{quotation}
[p. 6] \emph{\emph{``Theorem''.} There is a spectral sequence with
$E_2^{\Star,\Star} =$ ``something computable'' and converging to $H\Star$,
something desirable. The important observation to make about the statement of
the theorem is that it gives an $E_2$-term of the spectral sequence but says
nothing about the successive differentials $d_r$. Though $E_2^{\Star,\Star}$
may be known, without $d_r$ or some further structure, it may be impossible to
proceed.}

\ldots\ \ldots

[p. 28]  \emph{It is worth repeating the caveat about differentials mentioned
in Chapter 1: knowledge of $E_r^{\Star,\Star}$ and $d_r$ determines
$E_{r+1}^{\Star,\Star}$ but not $d_{r+1}$. If we think of a spectral sequence
as a black box, then the input is a differential bigraded module, usually
$E_1^{\Star,\Star}$, and, with each turn of the handle, the machine computes a
successive homology according to a sequence of differentials. If some
differential is unknown, then some other (any other) principle is needed to
proceed. From Chapter 1, the reader is acquainted with several algebraic tricks
that allow further calculation. In the non-trivial cases, it is often a deep
geometric idea that is caught up in the knowledge of a differential.}
\end{quotation}

It is in fact again a matter of computability. The higher differentials of a
spectral sequence are \emph{mathematically} defined, but, in most cases, their
definition \emph{is not} constructive: the differentials are not
\emph{computable} with the provided information. For example the result of
Adams' work about the first loop space is nothing but an algorithm computing
the higher differentials and solving the extension problems at abutment of the
corresponding Eilenberg-Moore spectral sequence, thanks to the coalgebra
structure over the initial cellular chain complex. But this does not compute
the coalgebra structure for the CW model of the loop space so that you cannot
continue: this is nothing but the ``severe'' difficulty above observed by
Carlsson and Milgram. See the nice work of Baues \cite{BAUS1} to go a little
further, but this does not give a solution for the general problem of
``iterating the cobar construction''.

\section{Three complete solutions for the computability problem.}\label{sec2}

In fact three solutions are now available to work in a \emph{constructive}
context in Algebraic Topology. This section describes the main ingredients of
the solutions that are due to Rolf Sch\"on \cite{SCHN} and Justin Smith
\cite{SMTH1,SMTH2}. The rest of the paper is devoted to the Rubio-S. solution
and the corresponding Kenzo program.

\subsection{Rolf Sch\"on's solution.}

Sch\"on's solution \cite{SCHN} is a systematic reorganization of Edgar Brown's
special work \cite{BRWNE1} around the computation of homotopy groups.
Frequently in Homological Algebra, we work with large chain complexes, the
homology groups of which are of finite type; for example the singular chain
complex of a compact manifold is not at all of finite type, but the homology
groups of this chain complex on the contrary are. The same in a simplicial
context; for example a simplicial \emph{group} version of the circle $S^1$
necessarily has an infinite number of simplices in any positive dimension, but
the homology groups are null or with only one generator. When you work with the
traditional tools of homological algebra, you must frequently handle highly
infinite chain complexes even if you know the final result is of finite type.

Edgar Brown designed an approximation process which has been skilfully
generalized by Rolf Sch\"on. Let $X$ be a simplicial set, described as the
limit of a sequence $(X_n)$ of finite approximations. Then the homology group
$H_p(X)$ is the inductive limit of the groups $(H_p(X_n))_n$, so that the
following definition could be useful.

\begin{definition}
--- A \emph{Sch\"on \Z-module} $G$ is a triple \[((G_n)_{n \geq 0}, (\phi_n)_{n
\geq 0}, \alpha)\] where the following conditions are satisfied. Every $G_n$ is
a \Z-module of finite type, and $\phi_n: G_n \rightarrow G_{n+1}$ is a morphism
of \Z-module; the sequence $(G_n, \phi_n:G_n \rightarrow G_{n+1})$ is an
inductive system and its limit~$G$ is again of finite type. The third component
$\alpha$ precisely describes how the limit is reached; $\alpha: \N \rightarrow
\N \times \N$ is as follows: if $\alpha(i) = (j,k)$, then $i \leq j \leq k$ and
the canonical morphism \(\mbox{Im}\,G_j \rightarrow G\) is in fact an
isomorphism:
\[
  \alpha: i \mapsto
    \left\{
      \begin{array}{ccl}
        G_j & \rightarrow & \mbox{Im}\,G_j \subset G_k \\
            &             & \mbox{\scriptsize$\cong$} \downarrow \\
            &             & \hspace{7pt}G
      \end{array}
    \right\}
\]
\end{definition}

The \emph{existence} of such a map $\alpha$ is implied by the finiteness
property of the inductive limit~$G$ which is assumed, but an \emph{effective}
knowledge of this map is required. Because you do not know a priori what
approximations $X_n$ of $X$ will be later required for some calculation, the
value $\alpha(i)$ must be \emph{computable} for any $i$. We call the map
$\alpha$ the \emph{convergence descriptor}.

The books of Homological Algebra are full of theorems of this sort:

\begin{theorem}
--- There is an exact sequence:
\[
\cdots \rightarrow A \stackrel{f}{\rightarrow} B \rightarrow C \rightarrow D
\stackrel{g}{\rightarrow} E \rightarrow \cdots
\]
\end{theorem}

The underlying idea is that if you know the \Z-modules $A$, $B$, $D$ and $E$,
then you should be able to guess the unknown module $C$. Of course you must in
fact also know the maps $f:A \rightarrow B$ and $g: D \rightarrow E$ to
determine the modules $\mbox{Coker}(f)$ and $\mbox{Ker}(g)$, giving a simpler
exact sequence:
\[
0 \rightarrow \mbox{Coker}(f) \rightarrow C \rightarrow \mbox{Ker}(g)
\rightarrow 0,
\]
and now you could have an extension problem in front of you, about which the
exact sequence says nothing at all! The situation is analogous with the
spectral sequences but usually much more complicated. It was exactly the
problem encountered by Serre when he was looking for the group $\pi_6 S^2$: the
unknown group was in an exact sequence at the end of a spectral sequence
between two groups $\Z_2$ and $\Z_6$, and a new idea is necessary to terminate.

On the contrary such a problem is entirely solved in the framework designed by
Rolf Sch\"on. The situation is now the following: the modules $A$, $B$, $D$ and
$E$ are four known \emph{Sch\"on modules}; the map $f$ is in fact a morphism of
inductive systems and in particular for every $n$ a morphism $f_n$ is defined
satisfying the usual properties; the same between the other components of the
exact sequence. For the \emph{unknown} Sch\"on module $C$, the underlying
inductive system is known but its convergence descriptor \emph{is not}. You
know there is an \emph{exact} sequence between the limits $A$, \(B\), \(C\),
\(D\) and $E$, but at the $n$-th stage of the inductive systems, you have only
a ``differential'' sequence:
\[
A_n \stackrel{f_n}{\rightarrow} B_n \stackrel{f'_n}{\rightarrow} C_n
\stackrel{g'_n}{\rightarrow} D_n \stackrel{g_n}{\rightarrow} E_n
\]
where two successive maps have a null composition, but this sequence is not
necessarily exact.
\[
\begin{array}{l}
  (A_n, \phi_n),\ \alpha_A : i \mapsto (j,k) \\
  \hspace{6pt} \downarrow \mbox{\scriptsize ($f_n$)} \\
  (B_n, \phi'_n),\ \alpha_B : i \mapsto(j,k) \\
  \hspace{6pt} \downarrow \mbox{\scriptsize ($f'_n$)} \\
  (C_n, \chi_n),\ \ ???????????? \\
  \hspace{6pt} \downarrow \mbox{\scriptsize ($g'_n$)} \\
  (D_n, \psi'_n),\ \alpha_D : i \mapsto (j,k) \\
  \hspace{6pt} \downarrow \mbox{\scriptsize ($g_n$)} \\
  (E_n, \psi_n),\ \alpha_E : i \mapsto (j,k)
\end{array}
\]

\begin{theorem}
\emph{(Sch\"on \cite{SCHN})} --- With the previous data, an algorithm can
compute the convergence descriptor of the intermediate Sch\"on module $C$.
\end{theorem}

Once the missing descriptor $\alpha_C$ is available, then you can compute the
limit $C$. But, and maybe this is more important, the process is \emph{stable}:
the object $C = ((C_n),(\chi_n), \alpha_C)$ which is returned by Sch\"on's
algorithm is again a Sch\"on module and can be a part of the input for another
call of the same algorithm. Rolf Sch\"on explains in his nice paper \cite{SCHN}
how this method allows to entirely transform classical Homological Algebra into
a \emph{constructive} theory.

To our knowledge, Sch\"on's work has not yet led to concrete machine programs.
It is a pity: his general framework is quite original and interesting with
respect to what is usually done in computational algebra. The opinion of the
present author is that concrete implementation of Sch\"on's results must
absolutely be done and should give new insights into several fields: at least
in symbolic computation, in computational algebra and also in algebraic
topology.

\subsection{Justin Smith' solution.}

This second solution is quite different from the previous one. In a sense it is
exactly the solution of the problem stated by Karoubi (cf Section~\ref{sec1}).
Let $X$ be a \emph{simplicial} complex. The main problem in Algebraic Topology
comes from the non-commutativity of the Alexander-Whitney diagonal. If you
intend to send an interval $I$ onto the diagonal of a square $I \times I$,
using only the \emph{bi}simplicial structure of this square, that is, using
\emph{only} its four boundary edges, then you can join one vertex to the
opposite one turning around the square in two different ways:
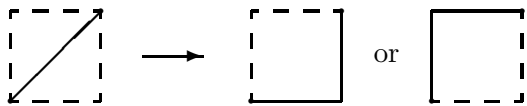
\setlength{\unitlength}{0.8cm}
\begin{center}
\begin{picture}(0,2)(0,-0.25)
  \thicklines
  \put(-4,0){\line(1,1){1.5}}
  \put(-4,0){\circle*{0.1}}
  \put(-2.5,1.5){\circle*{0.1}}
  \put(-1.8,0.75){\vector(1,0){1}}
  \put(-4,0){\dashbox{.25}(1.5,1.5){}}
  \put(0,0){\dashbox{.25}(1.5,1.5){}}
  \put(0,0){\line(1,0){1.5}}
  \put(1.5,0){\line(0,1){1.5}}
  \put(0,0){\circle*{0.1}}
  \put(1.5,1.5){\circle*{0.1}}
  \put(2.25,0.75){\makebox(0,0){or}}
  \put(3,0){\dashbox{.25}(1.5,1.5){}}
  \put(3,0){\line(0,1){1.5}}
  \put(3,1.5){\line(1,0){1.5}}
  \put(3,0){\circle*{0.1}}
  \put(4.5,1.5){\circle*{0.1}}
\end{picture}
\end{center}

These paths are different but they are homotopic. This homotopy is quite
important and leads to this diagram:

\setlength{\unitlength}{1pt}
\begin{center}
\begin{picture}(0,48)(0,0)
 \put(-60,+40){\makebox(0,0){\(C_\Star(X^2)\)}}
 \put(+60,+40){\makebox(0,0){\(C_\Star(X) \otimes C_\Star(X)\)}}
 \put(-60,0){\makebox(0,0){\(C_\Star(X^2)\)}}
 \put(+60,0){\makebox(0,0){\(C_\Star(X) \otimes C_\Star(X)\)}}
 \put(-35,+40){\vector(1,0){50}}
 \put(-35,0){\vector(1,0){50}}
 \put(-60,+30){\vector(0,-1){20}}
 \put(+60,+30){\vector(0,-1){20}}
 \put(-35,+28){\vector(3,-1){50}}
 \put(-5,+42){\makebox(0,0)[b]{\(\Delta\)}}
 \put(-5,+2){\makebox(0,0)[b]{\(\Delta\)}}
 \put(-62,+20){\makebox(0,0)[r]{\(\pi\)}}
 \put(+58,+20){\makebox(0,0)[r]{\(\pi\)}}
 \put(-4,+21){\makebox(0,0)[bl]{\(h\)}}
\end{picture}
\end{center}

The chain complex $C_\Star(X^2)$ is obtained from the canonical
\emph{simplicial} structure of $X^2$; on the contrary the other chain complex
$C_\Star(X) \otimes C_\Star(X)$ comes from the canonical \emph{bisimplicial}
structure of the same space. If for example $X$ is the interval $I = [0,1]$, in
the first case a square is presented as the union of two triangles joined along
a diagonal; in the second case no diagonal in the square, only the boundary
edges, the square is simply the product of two intervals. Both presentations
are related by the Alexander-Whitney map $\Delta$. Furthermore both components
of $X^2$ can be swapped, and this leads to the vertical canonical (different)
maps $\pi$. Then the diagram is not commutative: $\Delta \circ \pi \neq \pi
\circ \Delta$. Nevertheless the homotopy operator $h$ explains both maps are
homotopic. But the same difficulty occurs now for the homotopy $h$ which in
turn is not compatible with the symmetry of its source and its target, but
again a homotopy can be constructed and so on. This process roughly explained
here for both factors works also for an arbitrary number of factors $X^n$ and
all the homotopies are related by a very rich structure called a
\emph{coalgebra structure with respect to the symmetric operad $\mathfrak{S}$}.

Using an appropriate modified model for the symmetric operad $\mathfrak{S}$ and
also a corresponding notion of coalgebra called $m$-structure, Justin Smith
succeeded firstly in iterating the cobar construction \cite{SMTH1}, and more
recently \cite{SMTH2} in proving that a chain complex carrying an $m$-structure
contains a homotopy type, so that such a structure can be used as the
$\mathcal{H}$~component (cf Definition~\ref{computable-model}) for the
computable algebraic model demanded by Karoubi.

While preparing this paper, the author received a message of Justin Smith
announcing a partial programming work was just starting around the symmetric
operad \(\mathfrak{S}\). So that we can hope Justin Smith' solution finally
leads also to a concrete computer program. The situation here is also
interesting because of the original environment where work is to be undertaken:
it is probably the first time an operad structure is implemented. Certainly, at
least because they solve the same problem (!), Justin Smith' program and ours
will be strongly related. Probably the structure of Justin Smith's solution is
richer than for our solution; the latter works essentially like a blackbox,
because of its highly functional process which in a sense hides what actually
happens during the execution. When both solutions will be available,
determining what exactly the relations between them are will be still more
interesting!

\subsection{A quick sketch of the Rubio-S. solution.}

The nature of this ``third''\footnote{The first announcement goes back to 1987
\cite{SRGR1}; the first computer program computing an iterate cobar
construction started in 1990 \cite{RBSS1}.} solution is not so far from Justin
Smith' one. In our framework, any reasonable homotopy type is described as
follows: firstly a free \Z-chain complex of finite type in any dimension
\(EC_\Star X\) is given; then a further structure \(\mathcal{H}\) is added to
this chain complex in such a way a homotopy type is finally so defined; in fact
this homotopy type can be realized as a CW-complex, the cellular complex of
which being \(EC_\Star X\); it is well known this cellular complex does not
define a homotopy type, but the added structure \(\mathcal{H}\) gives the
missing information. What is quite original with respect to the traditional
organization in Algebraic Topology is the deeply \emph{functional} nature of
the structure \(\mathcal{H}\), the main subject of the rest of this paper.

\section{A didactic example of functional programming.}\label{sec3}

We briefly recall in this section a typical situation where it is much better
to work with functional objects carrying an enormous information, instead of
working with data close to those that are looked for.

Let \(G\) be a finite graph \(G = (V,E)\); the set \(V\) is the vertex set and
\(E\) is the set of the edges. A \emph{good} colouring of \(G\) consists in
defining a colour for each vertex so that two adjacent vertices have different
colours. The \emph{chromatic number} \(\chi(G)\) is the minimal number of
colours that are necessary. It is not so easy to design a program computing
this chromatic number. The traditional backtracking methods work but are quite
inefficient.

If you think of a recursive method, you cannot design such a method if you work
only with the chromatic number. Let \(\alpha \in E\) be an edge between the
vertices \(v,w \in V\). You would like for example to deduce \(\chi(G)\) from
\(\chi(G')\) where \(G'\) is the graph \(G\) without the edge~\(\alpha\). In
fact two interpretations of \(G'\) make sense. The first one \(G_1\) has the
same vertex set as~\(G\), and \(\alpha\) is simply removed from \(E\). The
second interpretation \(G_2\) consists in collapsing the edge \(\alpha\) over
one vertex coming for both vertices \(v\) and \(w\); in particular if we
previously had two different edges \(uv\) and \(uw\) starting from another
vertex \(u\) and going respectively to \(v\) and \(w\), both edges give only
one edge in \(G_2\): both \(G\)-edges are now identified in \(G_2\). For
example if \(G\) is a complete graph of order \(n\), then \(G_1\) is the same
with only one edge removed, but \(G_2\) is the complete graph of order \(n-1\).
And very simple cases show the knowledge of \(\chi(G_1)\) and \(\chi(G_2)\) is
not sufficient to compute \(\chi(G)\): the chromatic number does not contain
enough information; \emph{we need more}. \setlength{\unitlength}{1cm}
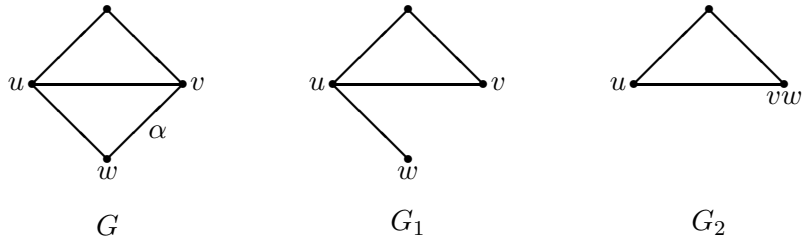
\begin{center}
\begin{picture}(0,3)(0,-2)
  \thicklines
  \put(-4,-2){\makebox(0,0)[b]{\(G\)}}
  \put(-5.1,0){\makebox(0,0)[r]{\(u\)}}
  \put(-2.9,0){\makebox(0,0)[l]{\(v\)}}
  \put(-4,-1.1){\makebox(0,0)[t]{\(w\)}}
  \put(-3.45,-0.55){\makebox(0,0)[lt]{\(\alpha\)}}
  \put(-4,+1){\circle*{0.1}}
  \put(-5,0){\circle*{0.1}}
  \put(-3,0){\circle*{0.1}}
  \put(-4,-1){\circle*{0.1}}
  \put(-5,0){\line(1,1){1}}
  \put(-5,0){\line(1,0){2}}
  \put(-5,0){\line(1,-1){1}}
  \put(-4,+1){\line(1,-1){1}}
  \put(-4,-1){\line(1,1){1}}
  \put(0,-2){\makebox(0,0)[b]{\(G_1\)}}
  \put(-1.1,0){\makebox(0,0)[r]{\(u\)}}
  \put(+1.1,0){\makebox(0,0)[l]{\(v\)}}
  \put(0,-1.1){\makebox(0,0)[t]{\(w\)}}
  \put(0,+1){\circle*{0.1}}
  \put(-1,0){\circle*{0.1}}
  \put(+1,0){\circle*{0.1}}
  \put(0,-1){\circle*{0.1}}
  \put(-1,0){\line(1,1){1}}
  \put(-1,0){\line(1,0){2}}
  \put(-1,0){\line(1,-1){1}}
  \put(0,+1){\line(1,-1){1}}
  \put(4,-2){\makebox(0,0)[b]{\(G_2\)}}
  \put(2.9,0){\makebox(0,0)[r]{\(u\)}}
  \put(5,-0.1){\makebox(0,0)[t]{\(vw\)}}
  \put(4,+1){\circle*{0.1}}
  \put(3,0){\circle*{0.1}}
  \put(5,0){\circle*{0.1}}
  \put(3,0){\line(1,1){1}}
  \put(3,0){\line(1,0){2}}
  \put(4,+1){\line(1,-1){1}}
\end{picture}
\end{center}
Let us consider the \emph{chromatic polynomial} \(P_G(X)\); it is a polynomial
with one variable defined as follows: if \(n\) is a positive integer, then
\(P_G(n)\) is the number of good colourings of \(G\) that are possible with
\(n\) colours. Now the situation is good:

\begin{quotation}
1) A recursive relation holds and it is simple: \(P_G(X) = P_{G_1}(X) -
P_{G_2}(X)\); in fact, let us consider a good colouring of \(G_1\); then,
depending on whether both colours of \(v\) and \(w\) are the same or not, you
obtain a good colouring for \(G_2\) or \(G\), and the relation between \(P_G\),
\(P_{G_1}\) and \(P_{G_2}\) follows ; starting with graphs without any edge,
you obtain in particular that \(P_G\) actually is a polynomial!

2) The polynomial \(P_G\) contains an \emph{infinite} number of elementary
data: how many good colourings exist with 1 color, with 2 colours, and so on;
we now have enough information;

3) These data are coded in a \emph{functional} way: of course you cannot store
in your machine all the values \(P_G(n)\); but it is sufficient to store the
degree and the coefficients of \(P_G\): a polynomial is a \emph{finite} object
which is nothing but a \emph{program} ready to compute the value \(P_G(n)\) for
every integer \(n\) in the \emph{infinite} set~\N.

4) The chromatic number \(\chi(G)\) is a \emph{by-product} of the polynomial
\(P_G\): it is sufficient to compute \(P_G(1)\), \(P_G(2)\), \ldots, until you
find the first integer \(n\) satisfying \(P_G(n) > 0\); then \(\chi(G) = n\).
\end{quotation}

It is then easy to write down a recursive program computing the chromatic
number; it is more efficient than a program using backtracking, but however it
has an exponential complexity; the problem of finding a polynomial time
algorithm computing the chromatic number is open: it is a special case of the
general \(NP\)-complete problem.

The Rubio-S. solution for constructive Algebraic Topology is quite similar. The
role of the chromatic number \(\chi(G)\) is played by an \emph{effective} chain
complex \(EC_\Star X\), which is the cellular chain complex of some CW-model of
the homotopy type we intend to algebraically define. The situation is the same:
the information given in this chain complex is in general too poor to process
new objects deduced from this one and others; \emph{we need more}. We will
define new ingredients, in general containing an \emph{infinite} number of
elementary data and which completely define a homotopy type; but these
ingredients will be coded in a \emph{functional} way so that a machine program
will be able to handle them as easily as polynomials\footnote{At least if your
programming language allows you to use \emph{functional programming}.} and to
compute the corresponding ingredients for a new homotopy type constructed from
others which were defined by means of such data.

\section{Objects with effective homology.}\label{sec4}

\subsection{Effective chain complexes.}

A chain complex is a sequence of \Z-modules and homomorphisms:
\[
\ldots \leftarrow C_{n-1} \leftarrow C_n \leftarrow C_{n+1} \leftarrow \ldots
\]
where the composition of two successive arrows is null.

\begin{fromnowon}
--- All the chain groups~\(C_n\) of a chain complex \(C_\Star\) are \emph{free
\Z-modules with distinguished basis}.
\end{fromnowon}

In the following definition, the set \(\mathcal{U}\) is the ``machine
universe'': any machine object is an element of \(\mathcal{U}\); the set
\(\mbox{List} \subset \mathcal{U}\) is the subset of all \emph{lists}, in other
words the finite sequences of elements of \(\mathcal{U}\).

\begin{definition}
--- An \emph{effective chain complex} is defined as a pair of algorithms:
\begin{itemize}
\item
\(\beta: \Z \rightarrow \mbox{List}\);
\iten
\(d: \Z \,\widetilde{\times}\,\mathcal{U} \rightarrow \mbox{List}\);
\end{itemize}
where:
\begin{enumerate}
\item
The output \(\beta(n)\) is the given basis of the free \Z-module \(C_n\); this
basis is a list and in particular is finite;
\iten
A pair \((n,g)\) is in \(\Z \,\widetilde{\times}\,\mathcal{U}\) if \(g\) is a
generator of \(C_n\), that is, if \(g \in \beta(n)\);
\iten
The output \(d(n,g)\) is a list representing the differential \(d_n(g) \in
C_{n-1}\).
\end{enumerate}
\end{definition}

If an effective chain complex \(C_\Star\) and an integer \(n\) are given, a
program can compute the boundary matrices in dimensions \(n+1\) and \(n\), and
an elementary algorithm then determines the homology group \(H_n(C_\Star)\).
The \emph{global} nature of an effective chain complex \(C_\Star\) is reachable
for any dimension \(n\).

\subsection{Locally effective chain complexes.}

\begin{definition}
A \emph{locally effective chain complex} \(C_\Star\) is defined as a pair of
algorithms:
\begin{itemize}
\item
\(\beta': \Z \times \mathcal{U} \rightarrow \mbox{Boolean}\);
\iten
\(d: \Z \,\widetilde{\times}\, \mathcal{U} \rightarrow \mbox{List}\);
\end{itemize}
where:
\begin{enumerate}
\item
The output \(\beta'(n,\gamma)\) is the Boolean \texttt{true} if and only if the
object \(\gamma\) is a generator of the chain group \(C_n\);
\iten
The sub-product \(\Z \,\widetilde{\times}\, \mathcal{U}\) is interpreted as in
the previous definition and the differential~\(d\) as well.
\end{enumerate}
\end{definition}

It is explained in the handbooks of set theory there are two different methods
to define a set \(S\). You can give the element list of \(S\); in a
computational framework, such a list is necessarily finite. You can also define
the set \(S\) by means of a characteristic property of its elements. For
example you can require an element of \(S\) must be an integer and must be odd.
Then such a set may be infinite. Do not object the set of actual elements that
can actually be processed on your machine is finite; consider for example this
Lisp definition of the set~\(\N_{\mbox{\scriptsize odd}}\):

\begin{lisp}
\begin{verbatim}
> (setf odd-integers
     #'(lambda (object)
          (and (integerp object)
               (oddp object))))
\end{verbatim}
\end{lisp}
\noindent This string of 82 characters is finite and \emph{defines} the
infinite set of odd integers.

In the same way the generators of our locally effective chain complexes are
defined by means of a characteristic property, so that now our chain groups are
not necessarily of finite type. This looks like an advantage with respect to
the notion of effective chain complex, but there is an important drawback: in
general no global information is reachable for such a machine chain complex; in
particular the homology groups in general are not computable. This is an avatar
of the main incompleteness theorem (G\"odel, Church, Turing, Post). The key
point of the Rubio-S. solution for Constructive Algebraic Topology consists in
combining effective \emph{and} locally effective chain complexes, connecting
them by \emph{reductions}.

\subsection{Reductions.}

\begin{definition}
A \emph{reduction} \(\rho: D_\Star \Rightarrow C_\Star\) between two chain
complexes is a triple \(\rho = (f,g,h)\) where:
\begin{enumerate}
\item
The components \(f\) and \(g\) are chain complex morphisms \(f: D_\Star
\rightarrow C_\Star\) and \(g: C_\Star \rightarrow D_\Star\);
\iten
The component \(h\) is a homotopy operator \(h: D_\Star \rightarrow D_\Star\)
(degree 1);
\iten The following relations are satisfied:
  \begin{enumerate}
  \item
  \(f \circ g = \mbox{id}_{C_{\scriptstyle *}}\);
  \(g \circ f + d_{D_{\scriptstyle *}} \circ h + h \circ d_{D_{\scriptstyle *}}
        = \mbox{id}_{D_{\scriptstyle *}}\);
  \item
  \(f \circ h = 0\); \(h \circ g = 0\) ; \(h \circ h = 0\).
  \end{enumerate}
\end{enumerate}
\end{definition}

In these formulas, \(d_{D_{\scriptstyle *}}\) denotes the differential of the
chain complex \(D_\Star\). These formulas have a simple interpretation: the
chain complex \(C_\Star\), the small one, is isomorphic to a subcomplex of
\(D_\Star\), the big one, and a decomposition \(D_\Star = C_\Star \oplus
E_\Star\) is given where the summand \(E_\Star\) is acyclic and provided with
an explicit homological contraction. This implies both chain complexes
\(C_\Star\) and \(D_\Star\) have the same homology.

Frequently in the Rubio-S. context, the big chain complex \(D_\Star\) is
locally effective, so that its homology groups are not computable; on the
contrary, the small chain complex \(C_\Star\) is effective, so that its
homology groups are computable. In such a situation, the reduction can be
understood as a provided description of the global homological properties of
\(D_\Star\). In particular if you are interested by the explicit value of
\(H_n(D_\Star)\), you can obtain the result by \(H_n(C_\Star)\); furthermore an
explicit representative for any homology class can be deduced in \(D_n\); if
\(z\) is a cycle of \(D_n\), the homology class of \(z\) can be determined, and
if null, a chain \(c \in D_{n+1}\) can be found such that \(dc = z\). In a word
you know \emph{everything} about the homological properties of \(D_\Star\).

\begin{definition}
An \emph{equivalence} \(\varepsilon: C_\Star \Leftarrow \Rightarrow E_\Star\)
is a pair \(\varepsilon = (\rho_\ell, \rho_r)\) of reductions \(\rho_\ell:
D_\Star \Rightarrow C_\Star\) and \(\rho_r: D_\Star \Rightarrow E_\Star\).
\end{definition}

Again, frequently the chain complexes \(C_\Star\) and \(D_\Star\) are only
locally effective and the third one \(E_\Star\) is effective; so that the
equivalence \(\varepsilon\) describes the homological properties of \(C_\Star\)
thanks to \(E_\Star\).

\subsection{Objects with effective homology.}

\begin{definition}
An \emph{object with effective homology} is a pair \((X,\varepsilon)\) where
\(X\) is some locally effective object and \(\varepsilon\) is an equivalence
between the chain complex ``canonically'' associated to \(X\) and some
effective chain complex.
\end{definition}

The associated chain complex depends on the context. For example if \(X\) is a
simplicial set, then \(C_\Star(X)\) could be the normalized chain complex
defining its simplicial homology. The simplicial set \(X\) should be also
locally effective; in other words some algorithm is given as the characteristic
property of the \(n\)-simplices of \(X\); if \(\sigma\) is such a simplex,
another algorithm can compute the faces \(\partial_i(\sigma)\). The
equivalence:
\[
\varepsilon: C_\Star(X) \stackrel{\rho_\ell}{\Leftarrow} D_\Star X
\stackrel{\rho_r}{\Rightarrow} E_\Star X
\]
entirely describes the homological properties of \(X\), because the chain
complex \(E_\Star X\) is effective. In general there is no way to \emph{deduce}
this equivalence from the locally effective object \(X\). Most often we start
with effective objects where such an equivalence is trivial, and also with
special objects for which the particular situation give such an equivalence;
the Eilenberg-MacLane spaces \(K(\pi,1)\) are of this sort if the group \(\pi\)
is abelian of finite type. Then the \emph{effective homology version} of the
``classical'' construction methods of Algebraic Topology allow you to obtain
new objects with effective homology. For example the Eilenberg-MacLane space
\(K(\pi,2)\) is the classifying space of \(K(\pi,1)\), so that the effective
homology version of the classifying space construction, available in the
program Kenzo, will give you a copy of \(K(\pi,2)\) with effective homology.
You can trivially iterate the process and obtain versions with effective
homology of the Eilenberg-MacLane spaces \(K(\pi,n)\)'s. Proceeding in the same
way with the loop space construction, a very simple solution for
\emph{iterating the cobar construction} is obtained.

\section{The spectral sequences revisited.}\label{sec5}

Many constructions in algebraic topology can be organized as solutions of
fibration problems. In particular the classifying space \(BG\) of a topological
group \(G\) is the solution for a fibration \(BG \times_\tau G\) where the
fiber space is the given group \(G\), the base space is the classifying space
\(BG\) and the product \(BG \times G\) is twisted in such a way the total space
\(BG \times_\tau G\) is contractible. The same idea where the base space \(X\)
is given and the fibre space is unknown leads to the loop space \(\Omega X\)
and the contractible total space \(X \times_\tau \Omega X\). The handbooks of
Algebraic Topology more or less explain the Eilenberg-Moore spectral sequence
can be used to ``compute'' the homology groups of the new objects \(BG\) and
\(\Omega X\) if the homology groups of \(G\) or \(X\) are known. In fact this
spectral sequence is in general unable to give you the new homology groups,
unless you are in a very special situation.

The Serre spectral sequence works in the third situation, when you are looking
for the homology groups of a total space \(B \times_\tau F\) if the homology
groups of \(B\) and \(F\) are known; but in general you meet the same
difficulties with the higher differentials and the extension problems at
abutment.

The Serre and Eilenberg-Moore spectral sequences have \emph{effective homology}
versions which work when the data are simplicial sets with effective homology.
We detail a little the organization and the proof for the Serre spectral
sequence.

\begin{theorem}
There exists an algorithm:
\begin{itemize}
\item
\emph{\textbf{Input:}} \parbox[t]{367pt}{Two simplicial sets \(B\) and \(F\)
with effective homology and a twisting operator~\(\tau\) defining a fibration
\(F \rightarrow B \times_\tau F \rightarrow B\);}
\iten
\emph{\textbf{Output:}} A version \emph{with effective homology} of the total
space \(T = B \times_\tau F\).
\end{itemize}
\end{theorem}

The same with the Eilenberg-Moore spectral sequences when you are looking for
the effective homology of the base space \(B\) (resp. the fiber space \(F\)),
if versions with effective homology of the total space \(T\) and the fiber
space \(F\) (resp. the base space \(B\)) are given. These effective homology
versions of the Serre and Eilenberg-Moore spectral sequences are available in
the program Kenzo.

The main ingredient for the proof of the effective homology version of the
Serre spectral sequence is the \emph{Basic Perturbation Lemma} \cite{BRWNR}.

\begin{theorem}
\emph{\textbf{Basic Perturbation Lemma}} --- Let \(\rho: D_\Star \Rightarrow
C_\Star\) be a chain complex reduction and \(\delta_{D_{\scriptstyle *}}:
D_\Star \rightarrow D_\Star\) a \emph{perturbation} of the differential
\(d_{D_{\scriptstyle *}}\) satisfying the nilpotency condition. Then a general
algorithm can compute a new reduction \(\rho': D'_\Star \Rightarrow C'_\Star\)
where the underlying graded modules of \(D_\Star\) and \(D'_\Star\) (resp.
\(C_\Star\) and \(C'_\Star\)) are the same, but the differentials are
perturbed:
\begin{eqnarray*}
  d_{D'_{\scriptstyle *}} &=&
           d_{D_{\scriptstyle *}} + \delta_{D_{\scriptstyle *}}\\
  d_{C'_{\scriptstyle *}} &=&
           d_{C_{\scriptstyle *}} + \delta_{C_{\scriptstyle *}}.
\end{eqnarray*}
\end{theorem}

The perturbation \(\delta_{D_{\scriptstyle *}}\) for the differential of the
big chain complex is \emph{given}; on the contrary the perturbation
\(\delta_{C_{\scriptstyle *}}\) for the small one is \emph{computed} by the
algorithm. In a sense, the perturbation of the big chain complex is also
\emph{reduced}. This is possible thanks to the nilpotency condition: let \(h:
D_\Star \rightarrow D_\Star\) be the homotopy component of the reduction
\(\rho\); then the nilpotency condition is satisfied if the composition \(\nu =
h \circ \delta_{D_{\scriptstyle *}}\) is pointwise nilpotent, that is,
\(\nu^n(x) = 0\) for an \(n \in \N\) depending on \(x\).

A typical application of the basic perturbation lemma is the following. Let \(T
= B \times_\tau F \) be a fibration with the the base space \(B\) and the fiber
space \(F\). Let us assume two reductions \(\rho_B: C_\Star(B) \Rightarrow
EB_\Star\) and \(\rho_F: C_\Star(F) \Rightarrow EF_\Star \) are given,
describing the homology of both spaces by means of the \emph{effective} chain
complexes \(EB_\Star\) and \(EF_\Star\); then it is easy, thanks to
Eilenberg-Zilber, to compute a \emph{non-twisted} product reduction:
\[
\rho_B \times \rho_F: C_\Star(B \times F) \Rightarrow EB_\Star \otimes
EF_\Star.
\]
The underlying graded modules of \(C_\Star(T) = C_\Star(B \times_\tau F)\) and
\(C_\Star(B \times F)\) are the same but the differentials are not; the
difference is a perturbation of the big chain complex. If the base space \(B\)
is 1-reduced (no edge, the geometry begins in dimension 2), then the nilpotency
condition is satisfied and applying the Basic Perturbation Lemma gives a
reduction:
\[
\rho_T : C_\Star(T) = C_\Star(B \times_\tau F) \Rightarrow EB_\Star \otimes_t
EF_\Star
\]
which describes the homology of the total space of the fibration by means of a
twisted tensor product of the chain complexes \(EB_\Star\) and \(EF_\Star\).

This was already done by Shih \cite{SHIH} and the present work about effective
homology is nothing but the following remark: if functional programming is
used, then Shih's presentation of the Serre spectral sequence becomes an
\emph{algorithm} computing a version \emph{with effective homology} of the
total space of a fibration if analogous versions of the fibre and base spaces
are given, at least if the base space is simply connected. It is a little more
complicated but not very difficult to process in the same way the
Eilenberg-Moore spectral sequences to compute a version with effective homology
of the base space or the fibre space if such versions of both other components
of the fibration are given.

\section{Computing homotopy groups.}\label{sec5.5}

\begin{theorem}
--- Let \(X\) be a 1-reduced (one vertex, no edge) simplicial set with effective
homology. Then the homotopy groups of \(X\) are computable.
\end{theorem}

This is a strong generalization of Edgar Brown's theorem about the
computability of homotopy groups of finite 1-reduced simplicial sets
\cite{BRWNE1}. Furthermore our proof is not difficult and leads to concrete
programs actually computing the first homotopy groups of a ``reasonable''
simplicial set; an example is given in Section \ref{sec8}.

Let \(\pi = \pi_n X\) the first non-zero homotopy group. Hurewicz' theorem
implies this group is also the first non-trivial homology group \(H_n(X,\Z) =
\pi \), a group which is computable, because \(X\) has effective homology. Then
a fundamental cohomology class \(\zeta \in H^n(X,\pi)\) is defined, which in
turn defines a canonical fibration:
\[
K(\pi, n-1) \hookrightarrow X_{n+1} \rightarrow X.
\]

The group \(\pi\) is of finite type so that starting from \(K(\pi, 1)\) and
using \((n-2)\) times the version with effective homology of the
Eilenberg-Moore spectral sequence gives a copy with effective homology of
\(K(\pi, n-1)\). Then applying our version of the Serre spectral sequence
produces the total space \(X_{n+1}\) of our fibration with its effective
homology. This total space is the same space as \(X\) except that the \(n\)-th
homotopy group is null: \(\pi_n X_{n+1} = 0\). Applying again Hurewicz' theorem
to \(X_{n+1}\) gives \(\pi_{n+1} X = \pi_{n+1} X_{n+1} = H_{n+1}(X_{n+1},
\Z)\). Iterating the process gives the result.

This sequential process to compute the homotopy groups is known as the
\emph{Whitehead tower}. The dual process (\emph{Postnikov tower}) may be used
as well, computing also the \emph{Postnikov invariants}.

\section{The Kenzo program.}\label{sec6}

The \emph{Kenzo} program implements the main components of the organization
that is roughly described in these notes. It is a 16000 lines Lisp program,
www-reachable at the address~\cite{DSSS}, with a rich documentation (340pp.).
It can be used with any Common Lisp system satisfying the ANSI
norm\footnote{Mainly Allegro Common Lisp (cf \texttt{www.franz.com}), LispWorks
(\texttt{www.harlequin.com}) and Mac Common Lisp (\texttt{www.digitool.com}).}.
A small typical demonstration is www-visible \cite{DSSS}.

It seems difficult to realize the same work with another programming language.
At least for four reasons:

\begin{itemize}
\iten The heart of our programming work is mainly devoted to complex functional
programming; this feature forbids to use the so called imperative languages
such as C++ or Java with which functional programming is theoretically
possible\footnote{All languages are ``equivalent''.}, but practically it is
not.
\iten The structures of Algebraic Topology that are processed by the Kenzo
program are rich and complex: chain complexes, differential graded algebras,
differential coalgebras, differential Hopf algebras, simplicial sets, Kan
simplicial sets, simplicial groups, various morphisms between these objects,
reductions, equivalences between chain complexes. In the current context, the
modern methods of Object Oriented Programming (OOP) \emph{must be used}. In
particular the multi-inheritance feature available in Common Lisp is
invaluable: for example a simplicial group is simultaneously a simplicial set
and a differential graded algebra, and these classes are both subclasses of the
class of chain complexes. In functional programming languages such as ML or
Maple-V\footnote{Functional programming is available in Maple-V release 5.},
the OOP tools that are provided are too weak (or lacking) to work comfortably.
On the contrary, from this point of view, Axiom would be satisfactory,
but\ldots
\iten The time complexity of the algorithms implemented in the Kenzo program is
high; more simply, computing time is critical. Common Lisp is a stratified
language where the lowest level can be understood as the assembly language of a
virtual machine (functions \texttt{car}, \texttt{cdr}, \texttt{cons},\ldots)
and the Lisp compiler produces very efficient code for the low level functions.
So that using this assembly-like language when programming the kernel of a
program is an excellent optimization tool. Furthermore the powerful Lisp
macrogenerator allows the user to define his own intermediate language. Other
good languages such as Axiom, ML, Maple have a too thick interface between the
machine and the user to be satisfactory from this point of view.
\iten Lisp is one of the oldest languages still available and his enormous
and well organized package of predefined functions, for example to process
lists, trees, binary numbers, gives the user powerful tools again not available
in the other current high level languages, in particular when dynamically
created functions are implied.
\end{itemize}

No particular difficulty has been met during the programming work. In
particular, the rigorous \emph{mathematical} definition of the virtual Common
Lisp machine \cite{STEL,CLHS} gives the programmer a safe and convenient
framework.

\section{New research fields.}\label{sec7}

Various \emph{new} research fields are open by this work, in computer science
and in ``pure'' mathematics as well. Let us quickly describe two typical
examples.

\subsection{A new subject in Computer Science.}

A Kenzo computation of some homology group, for example a homology group of an
iterated loop space \(H_p \Omega^n X\) is split in two steps :

\begin{enumerate}
\iten Constructing a version \emph{with effective homology} of the loop space \(\Omega^n
X\); during this step, an enormous set of functional objects, something like
several hundreds or thousands, are dynamically constructed. They are organized
as an oriented graph where the nodes are the functional objects and each node
\(f\) is connected to several other nodes \(f_1,\ldots,f_k\) if a call of \(f\)
requires the call of \(f_1,\ldots,f_k\), to be viewed as auxiliary functions
(subroutines), which in turn have other auxiliary functions, and so on. But at
this time these functions have not yet worked: the first step is in a sense
\emph{macrogeneration} of \emph{object}\footnote{In fact, this is an illusion:
thanks to the \emph{closure} mechanism, only an enormous set of pointers is
installed.} code;
\iten When the computation of \(H_p \Omega^n X\) is started, the effective
chain complex corresponding to \(\Omega^n X\) is examined, two (finite)
boundary matrices are constructed, and the homology group is computed. The
construction of this boundary matrix is the problem with ``severe''
difficulties mentioned by Carlsson and Milgram, see Section \ref{sec1}; the
``program''  written in the step 1 now works and most functions are used.
\end{enumerate}

This situation gives rise to a difficult and interesting problem of memory
optimization. When the function \(f\) is called and some result
\(f(x_1,\ldots,x_k)\) has been computed, what about the idea of storing the
result? After all, and this is frequent, the same calculation will be again
required later. If the calculation is trivial, for example if the map \(f\) is
constant, or if it is fast, storing the result is expensive in time and space.
If on the contrary the computation is long, it is better to store the result to
avoid the repetition. But the decisions that are to be taken are not
independent from each other: if the calculation of \(f(x)\) is long but amounts
in fact to calculating \(f_1(x')\), storing the result \(f_1(x')\) implies the
calculation of \(f(x)\) becomes very fast! Furthermore, after a long work,
\emph{experience} can show that in fact some stored result has never been
reused, so that it could be thrown away? Yes, but in general the program is
unable to prove the result will \emph{certainly} not be re-used. It seems clear
only empirical methods can be applied, but nevertheless modelizing and studying
simplified models from this point of view should be interesting and useful.

In the Kenzo program, a small set of empirical methods are applied to decide
when a result is stored or not, but it is obvious we are far from the ``best''
choices.

\subsection{A new research field in pure mathematics.}

The complicated calculations which may be undertaken with the help of the Kenzo
program give new insights into some fields. The following example is typical.
If \(X\) is a 1-reduced (one vertex, no edge) simplicial set, the main result
which was obtained by Adams~\cite{ADMS}\footnote{See also \cite{CRML} for an
excellent recent extensive study of the subject.} towards the calculation of
the homology groups \(H_\Star \Omega X\) was a morphism of differential graded
algebras:
\[
\alpha: \mbox{Cobar}^{C_{\scriptstyle \ast}X}(\Z,\Z)
        \longrightarrow
        C_\Star \Omega X
\]
which is a chain equivalence. In interesting cases, the source of \(\alpha\) is
of finite type. The computation of \(H_\Star \Omega X\) amounts to considering
the chain complex \(\mbox{Cobar}^{C_{\scriptstyle \ast}X}(\Z,\Z)\) and its
finiteness properties make the homology groups computable. The Kenzo program
computes such a map \(\alpha\) and also an \emph{explicit} inverse chain
equivalence:
\[
\beta:  C_\Star \Omega X
        \longrightarrow
        \mbox{Cobar}^{C_{\scriptstyle \ast}X}(\Z,\Z).
\]

Once upon a time, a student implicitly used that \(\beta\) is also a morphism
of differential graded algebra. To persuade him he was wrong, the present
author used  the Kenzo program to give him simple examples showing such a
statement is not sensible, but he was rather surprised: the map~\(\beta\)
automatically constructed by the Kenzo program is, at least for the numerous
examples that have been tried, a morphism of algebra! In fact so many cases
have been computed that this is now an \emph{experimental} ``definitive''
\emph{fact}. This is an amazing strong version of Adams' result: there exists a
two-sided ideal \(I\) in the algebra \(C_\Star \Omega X\) such that Adams'
Cobar construction \(\mbox{Cobar}^{C_{\scriptstyle \ast}X}(\Z,\Z)\) is nothing
but the quotient \(C_\Star \Omega X/I\).

This became the main research subject of this student. Several interesting
results in this direction have been obtained, but at this time, the complete
result has not yet been proved. In particular it was completely obtained if a
new differential is installed on \(\mbox{Cobar}^{C_{\scriptstyle
\ast}X}(\Z,\Z)\), but it is not clear what the status of this new differential
is. See \cite{DNCT1,DNCT2}.

Other amazing experimental results of this sort have been obtained, in
particular around the canonical \emph{algebraic} fibration:
\[
C_\Star \Omega X \rightarrow X \otimes_t C_\Star \Omega X
                 \rightarrow X.
\]

This is the \emph{algebraic} version of the co-universal fibration:
\[
\Omega X \hookrightarrow PX \rightarrow X
\]
where the fibre space (resp. total space) is the loop space (resp. the path
space) of the pointed space X. The path space is contractible: it is a ``unit''
space and in a sense, \(\Omega X\) is an inverse space of \(X\). In the same
way, the twisted tensor product \(X \otimes_t C_\Star \Omega X\) is acyclic and
an explicit contraction \(h\) of this chain complex plays a capital role in
effective homology. The \emph{existence} of this contraction is known for a
long time \cite{BRWNE2}, but the explicit Kenzo computation of \(h\) shows very
surprising properties, which imply we are far from mastering the underlying
algebraic structure. Let us recall the loop space construction is the heart of
Algebraic Topology and that many problems can be reduced to problems about loop
spaces; they were \emph{invented} by Jean-Pierre Serre fifty years ago for this
reason.

\section{Examples of calculations.}\label{sec8}

\subsection{\(H_5 \Omega^3 \mathtt{Moore}(\mathbb{Z}_2, 4)\).}

Carlsson and Milgram explain in the paper quoted in Section \ref{sec1} the
computation of \(H_\Star \Omega^n X\) may be undertaken if \(X\) is a
suspension \(X = S^n Y\); then the homology groups \(H_\Star \Omega^n X\) are
entirely determined by the homology groups \(H_\Star Y\) thanks to a process
where the Dyer-Lashof homology operations play the main role, see
\cite{CRML,CHLM}. For example the Moore space \(\mbox{Moore}(\Z_2,4)\) is
nothing but the third suspension \(S^3P^2\R\), so that the homology groups
\(H_\Star \Omega^3 \mbox{Moore}(\Z_2,4)\) are entirely determined by the well
known groups \(H_\Star P^2\R = (\Z, \Z_2, 0, 0, \ldots)\). The best specialists
have been questioned and so far they have not yet been able to compute for
example \(H_5 \Omega^3 \mbox{Moore}(\Z_2,4)\)\footnote{In a case, two different
(!) results were successively proposed but both were wrong\ldots}. With the
Kenzo program the Moore space \(\mbox{Moore}(\Z_2,4) = S^3 P^2\R\) is
constructed as follows:

\begin{lisp}
\begin{verbatim}
  USER(3): (setf moore-2-4 (moore 2 4))
  [K1 Simplicial-Set]
\end{verbatim}
\end{lisp}

The (sub-) statement \texttt{(moore 2 4)} constructs the Moore space and the
statement \mbox{\texttt{(setf ...)}} assigns the result to the symbol
\texttt{moore-2-4}. Lisp explains the result is the Kenzo object \#1
(\texttt{[K1 \ldots]}) and this object is a simplicial set. Then the third loop
space is constructed and the result is assigned to the symbol
\texttt{o3-moore-2-4}:

\begin{lisp}
\begin{verbatim}
  USER(4): (setf o3-moore-2-4 (loop-space moore-2-4 3))
  [K30 Simplicial-Group]
\end{verbatim}
\end{lisp}

This time, the result is a simplicial \emph{group}. And the group \(H_5
\Omega^3 X = \Z_2^5\) is obtained in one minute:

\begin{lisp}
\begin{verbatim}
  USER(5): (homology o3-moore-2-4 5)
  Computing boundary-matrix in dimension 5.
  Rank of the source-module : 23.
  ;; Clock -> 1999-08-10, 14h 19m 56s.
  [... ... Lines deleted ... ...]
  Computing boundary-matrix in dimension 6.
  Rank of the source-module : 53.
  [... ... Lines deleted ... ...]

  Homology in dimension 5 :
  Component Z/2Z
  Component Z/2Z
  Component Z/2Z
  Component Z/2Z
  Component Z/2Z
  ---done---
  ;; Clock -> 1999-08-10, 14h 20m 50s.
\end{verbatim}
\end{lisp}

The Kenzo program has constructed a chain equivalence between the highly
infinite chain complex \(C_\Star \Omega^3 X\) and an effective one \(EC_\Star\)
which for example has 53 generators in dimension 5. The boundary matrices can
be computed and the corresponding homology group is obtained.

\subsection{A \(CW\)-model for \(\Omega^3 (P^\infty \mathbb{R} / P^3
\mathbb{R})\).}

Let us now consider an example where the Kenzo program overcomes the ``severe
difficulties'' quoted by Carlsson and Milgram, see again Section \ref{sec1}. In
a sense, the first case where their proposed methods fail is the following:
what about a CW-model for \(\Omega^3 X\) where \(X\) is the quotient \(X =
P^\infty \R / P^3 \R\)? Let us construct such a model with the Kenzo program;
the space \(X\) is constructed as follows:

\begin{lisp}
\begin{verbatim}
  USER(6): (setf p4 (r-proj-space 4))
  [K405 Simplicial-Set]
\end{verbatim}
\end{lisp}

The statement \texttt{(r-proj-space 4)} constructs the infinite real projective
space ``beginning'' only in dimension 4, that is the required quotient \(X =
P^\infty \R / P^3 \R\). The third loop space is constructed as before:

\begin{lisp}
\begin{verbatim}
  USER(7): (setf o3p4 (loop-space p4 3))
  [K434 Simplicial-Group]
\end{verbatim}
\end{lisp}

The Kenzo object \texttt{o3p4} is a simplicial group with effective homology
and the \emph{effective} associated chain-complex can be extracted:

\begin{lisp}
\begin{verbatim}
  USER(8): (setf eff-chain-complex-of-o3p4 (echcm o3p4))
  [K794 Chain-Complex]
\end{verbatim}
\end{lisp}

You see \(794 - 434 - 1 = 359\) other Kenzo objects (chain complexes with
various added structures and chain complex morphisms) have also been
constructed to obtain the result. The boundary matrix in dimension 5 of this
effective chain complex is computed by the Kenzo program in 30 seconds:

\begin{lisp}
\begin{verbatim}
  USER(9): (chcm-mat eff-chain-complex-of-o3p4 5)
  Computing boundary-matrix in dimension 5.
  Rank of the source-module : 33.
  ;; Clock -> 1999-08-10, 14h 22m 30s.
  [... ... Lines deleted ... ...]
  ;; Clock -> 1999-08-10, 14h 22m 57s.

  ========== MATRIX 13 lines + 33 columns =====
  L1=[C1=-2]
  L2=[C1=-1]
  L3=[C1=-4][C2=1][C3=-1][C4=-2]
  L4=[C2=1][C3=-1][C6=2]
  L5=[C1=6][C4=1][C6=1]
  L6=[C1=4][C4=4][C6=4][C7=3]
  L7=[C1=4][C12=-2][C14=2]
  L8=[C1=6][C4=1][C6=1]
  L9=[C1=4][C4=4][C6=4][C7=3]
  L10=[C8=4][C10=1][C11=-1][C14=-4][C15=-2][C20=-2]
  L11=[C1=4][C8=4][C10=1][C11=-1][C16=-4][C18=-1][C19=1][C23=-2]
  L12=[C12=4][C13=2][C16=-4][C18=-1][C19=1][C27=-2]
  L13=[C1=-1][C20=4][C21=2][C23=-4][C24=-2][C27=4][C28=2]
  ========== END-MATRIX
\end{verbatim}
\end{lisp}

You must read the result as follows : the non-null \(a_{i,j}\) terms of the
matrix are \(a_{1,1} = -2\), \(a_{2,1} = -1\), \ldots, \(a_{13,28} = 2\). This
is a computer-aided proof that there exists a CW-model for \(\Omega^3 X\) with
in particular 13 4-cells and 33 5-cells. This is an easy consequence of Adams'
Cobar construction, but the severe difficulties about the differentials are
here solved. In particular the boundary of the first 56cell \(e^5_1\) is
\(de^5_1 = -2e^4_1 - e^4_2 -4e^4_3 + 6e^4_5 + 4e^4_6 + 4e^4_7 + 6e^4_8 + 4e^4_9
+ 4e^4_{11} - e^4_{13}\). This defines only the homology type of the attaching
map for \(e^5_1\), but the rest of the Kenzo object contains also its
\emph{homotopy} type.

\subsection{\(\pi_5(\Omega S^3 \cup_2 e^3)\).}

The Kenzo program may computes the first homotopy groups of an \emph{arbitrary}
simply connected simplicial set with effective homology. Our last example of
Kenzo computation shows the calculation of \(\pi_5(\Omega S^3 \cup_2 e^3)\): a
3-cell \(e^3\) is attached to the loop space \(\Omega S^3\) by a map \(\partial
e^3 = S^2 \rightarrow \Omega S^3\) of degree 2. The space \(X = \Omega S^3
\cup_2 e^3\), called \texttt{dos3} below, can be constructed by a process which
is not necessary to detail here and which finishes as follows:

\begin{lisp}
\begin{verbatim}
  USER(13): (setf dos3 (disk-pasting os3 3 'new faces))
  [K826 Simplicial-Set]
\end{verbatim}
\end{lisp}

In principle the group \(H_2 X\) should be \(Z_2\):

\begin{lisp}
\begin{verbatim}
  USER(14): (homology dos3 2)
  Computing boundary-matrix in dimension 2.
  [... ... Lines deleted ... ...]
  Homology in dimension 2 :
  Component Z/2Z
  ---done---
\end{verbatim}
\end{lisp}
\noindent and the notion of a canonical cohomology class in dimension 2 is
defined ; the Kenzo program can construct it:

\begin{lisp}
\begin{verbatim}
  USER(15): (setf ch2 (chml-clss dos3 2))
  [K947 Cohomology-Class (degree 2)]
\end{verbatim}
\end{lisp}
The canonical fibration \(K(\Z_2,1) \hookrightarrow X_3 \rightarrow X\) induced
by this cohomology class is then constructed, and the total space of the
fibration is extracted:

\begin{lisp}
\begin{verbatim}
  USER(16): (setf f2 (z2-whitehead dos3 ch2))
  [K962 Fibration]
  USER(17): (setf X3 (fibration-total f2))
  [K968 Simplicial-Set]
\end{verbatim}
\end{lisp}

This is the beginning of the classical Whitehead tower, see Section
\ref{sec5.5}. In particular the group \(H_3 X_3 = \pi_3 X_3 = \pi_3 X\) can be
computed; in fact the Kenzo program has applied the version with effective
homology of the Serre spectral sequence:

\begin{lisp}
\begin{verbatim}
  USER(18): (homology X3 3)
  Computing boundary-matrix in dimension 3
  [... ... Lines deleted ... ...]
  Homology in dimension 3 :
  Component Z/2Z
  ---done---
\end{verbatim}
\end{lisp}
\noindent so that \(\pi_3 X = \Z_2\). Continuing in the same way for the
following stages of the Whitehead tower, the groups \(\pi_4 X = \Z + \Z_4\),
\(\pi_5 X = \Z_2^4\) are obtained in less than one hour.

\end{document}